\documentclass[12pt,a4paper,reqno]{amsart}
\usepackage{amsmath}
\usepackage{amsfonts}
\usepackage{amssymb}
\usepackage{bm}

\newcommand\R{{\mathbf{R}}}
\newcommand\C{{\mathbf{C}}}

\renewcommand\Re{{\operatorname{Re}}}

\theoremstyle{plain}
  \newtheorem{theorem}[subsection]{Theorem}
  \newtheorem{conjecture}[subsection]{Conjecture}
  
  \newtheorem{lemma}[subsection]{Lemma}

\theoremstyle{remark}
  \newtheorem{remark}[subsection]{Remark}

\theoremstyle{definition}

\include{psfig}

\begin{document}

\title[Two remarks on gKdV]{Two remarks on the generalised Korteweg de-Vries equation}
\author{Terence Tao}
\address{Department of Mathematics, UCLA, Los Angeles CA 90095-1555}
\email{tao@math.ucla.edu}
\subjclass{35Q53}

\vspace{-0.3in}
\begin{abstract}  We make two observations concerning the generalised Korteweg de Vries
equation $u_t + u_{xxx} = \mu ( |u|^{p-1} u )_x$.  Firstly we give a scaling argument that shows, roughly speaking,
that any quantitative scattering result for $L^2$-critical equation ($p=5$) automatically implies an analogous scattering
result for the $L^2$-critical nonlinear Schr\"odinger equation $iu_t + u_{xx} = \mu |u|^4 u$.  Secondly, in the
defocusing case $\mu > 0$ we present a new dispersion estimate which asserts, roughly speaking, that energy moves to the left faster than the mass, and hence strongly localised soliton-like behaviour at a fixed scale cannot persist for arbitrarily long
times.
\end{abstract}

\maketitle

\section{Introduction}

Consider the generalised Korteweg de Vries (gKdV) equation
\begin{equation}\label{eq:gkdv}
\partial_t u + \partial_{xxx} u = \mu \partial_x ( |u|^{p-1} u ), 
\end{equation}
where $p > 1$, $\mu$ is a real number, and $u: \R \times \R \to \R$ is a real-valued function.  The case $\mu=+1$ is defocusing, while the case $\mu=-1$ is focusing.  Formally at least, solutions to this equation have a conserved 
\begin{equation}\label{mass}
M(u) := \int_\R u(t,x)^2\ dx
\end{equation}
and a conserved energy
\begin{equation}\label{energy}
E(u) := \int_\R \frac{1}{2} [\partial_x u(t,x)]^2 + \frac{\mu}{p+1} |u(t,x)|^{p+1}\ dx
\end{equation}
as well as a scaling symmetry
$$ u(t,x) \mapsto \frac{1}{\lambda^{2/(p-1)}} u(\frac{t}{\lambda^3}, \frac{x}{\lambda} )$$
for any $\lambda > 0$.  The \emph{$L^2$-critical case} $p=5$ is of special interest, as the mass then becomes invariant under the above scaling.

In the focusing case $\mu = -1$ there is the \emph{soliton solution} $u(x) := Q(x-t)$, where the \emph{ground state}
$$ Q(x) = Q_p(x) := \left(\frac{p+1}{2 \cosh^2(\frac{p-1}{2} x)} \right)^{1/(p-1)}$$
is the unique positive even Schwartz solution $Q: \R \to \R^+$ to the ODE
$$ \partial_{xx} Q + Q^p = Q.$$
Further soliton solutions can then be created using the scaling symmetry, as well as the spatial translation symmetry $u(t,x) \mapsto u(t,x-x_0)$.

In this paper we shall establish two (unrelated) results on these equations.  The first is a link between 
the $L^2$-critical gKdV equation
and the $L^2$-critical NLS equation, which roughly speaking asserts that the former is at least as hard to control than the latter.
The second is a monotonicity formula for arbitrary defocusing gKdV equations, which asserts that the energy moves to the left faster than the mass, and hence strongly localised soliton-like solutions cannot occur.

\subsection{Scattering conjectures}

An $L^2_x(\R)$ solution to the free Airy equation
\begin{equation}\label{airy}
 \partial_t u + \partial_{xxx} u = 0
 \end{equation}
(which is the case $\mu=0$ of \eqref{eq:gkdv})
is known to obey the mixed norm estimate
$$ \| u \|_{L^5_x L^{10}_t(\R \times \R)} \leq C M(u)^{1/2} $$
for some absolute constant $C$; see \cite{kpv} (and \eqref{kato} below).  It is natural to ask whether estimates of this form are also available
in the nonlinear setting $\mu \neq 0$, or more precisely whether one has an \emph{a priori} bound of the form
$$
 \| u \|_{L^5_x L^{10}_t(I \times \R)} \leq f( M(u) )
$$
for all compact time intervals $I$ and all smooth\footnote{When $p$ is not an odd integer, one will have to settle for solutions in a moderately rough class, such as the energy class $C^0_t H^1_x(I \times \R)$.  Fortunately, we have a satisfactory $H^1$ wellposedness theory (see \cite{kpv}) which allows one to manipulate energy class solutions as if they were smooth.  In any event we shall primarily be interested in the case $p=5$.} solutions $u: I \times \R \to \R$ to \eqref{eq:gkdv} of a given mass $M(u)$, where $f: \R^+ \to \R^+$ is some function.
Scaling considerations suggest that such a claim is only possible in the $L^2$-critical case $p=5$; the soliton solution above then also shows that when $\mu=+1$, the claim is only possible when $M(u) < M(Q)$, i.e. when the mass of the solution is strictly less than that of the ground state\footnote{Furthermore, when $M(u) > M(Q)$, it is known that solutions can blow up in finite time \cite{merle}.  When $M(u) < M(Q)$ and the initial data has finite energy then it is known that the solution exists globally \cite{kpv}, \cite{weinstein}.}.  More precisely, we have

\begin{conjecture}[Scattering for critical gKdV]\label{gkdv-conj} Let $p=5$ and $\mu = \pm 1$.  Then there exist functions $f_+: [0,+\infty) \to [0,+\infty)$ and $f_-: [0,M(Q)) \to [0,+\infty)$ such that we have the bound 
$$
 \| u \|_{L^5_x L^{10}_t(I \times \R)}  \leq f_\mu( M(u) )
$$
for all compact time intervals $I$ and all Schwartz\footnote{By this we mean that $u$ is smooth, and all the spacetime derivatives of $u(t)$ are Schwartz in space, (locally) uniformly in time.} solutions $u: I \times \R \to \R$ to \eqref{eq:gkdv}.
\end{conjecture}

This conjecture, if true, would imply very strong qualitative and global control on energy class solutions to \eqref{eq:gkdv} of arbitrary finite mass in the defocusing case, or mass less than that of the ground state, for instance that they not only exist globally in time, but also scatter in the energy norm $H^1_x(\R)$ to solutions of the Airy equation \eqref{airy} as $t \to +\infty$ or $t \to -\infty$.  This essentially follows from the well-posedness theory in \cite{kpv}, and we omit the standard details.  This theory already gives the conjecture in the small mass case $M(u) \leq \epsilon_0$ for some small absolute 
constant $\epsilon_0 > 0$, but the conjecture is not known to be true in general.

The first result of this paper is to demonstrate the ``hardness'' of Conjecture \ref{gkdv-conj} by linking it to the analogous conjecture for the $L^2$-critical nonlinear Schr\"odinger (NLS) equation\footnote{The minus sign in front of the $i$ is convenient for us here, but one can easily reverse time (or conjugate $u$) to eliminate this sign if required.}
\begin{equation}\label{nls}
 -i\partial_t u + \partial_{xx} u = \mu |u|^{p-1} u
\end{equation}
where now $u: \R \times \R \to \C$ is complex-valued, and again $p=5$.  The mass \eqref{mass} and energy \eqref{energy} are also conserved quantities for this equation.  Note that in the focusing case $\mu=-1$ one has the soliton solution
$u(t,x) := e^{it} Q(x)$, which has the same mass $M(Q)$ as the soliton for gKdV.  $L^2$ solutions to the free Schr\"odinger
equation
$$ i\partial_t u + \partial_{xx} u = 0$$
are known to obey the Strichartz estimate
$$ \| u \|_{L^6_{t,x}(\R \times \R)} \leq C M(u)^{1/2};$$
see \cite{strich}.  Thus we are led to the counterpart for Conjecture \ref{gkdv-conj} for NLS:

\begin{conjecture}[Scattering for critical NLS]\label{nls-conj} Let $p=5$ and $\mu = \pm 1$.  Then there exist functions $f_+: [0,+\infty) \to [0,+\infty)$ and $f_-: [0,M(Q)) \to [0,+\infty)$ such that we have the bound 
$$
 \| u \|_{L^6_{t,x}(I \times \R)} \leq f_\mu( M(u) )
$$
for all compact time intervals $I$ and all Schwartz solutions $u: I \times \R \to \R$ to \eqref{nls}.
\end{conjecture}

Like Conjecture \ref{gkdv-conj}, Conjecture \ref{nls-conj} would have strong consequences for the global behaviour of this flow, but remains totally open (except when the mass is smaller than an absolute constant, see \cite{cwI}).  However, the four-dimensional analogue of this conjecture has recently been established in the radial case \cite{tao-visan}, while the analogue of this question for the energy-critical power $p = 1 + \frac{4}{d-2}$ (with $d \geq 3$) was settled in \cite{borg:book}, \cite{gopher}, \cite{rv}, \cite{visan}.  Thus there appears to be more progress on Conjecture \ref{nls-conj} than on Conjecture \ref{gkdv-conj}
at present.  Note also that unlike the gKdV, the NLS enjoys an additional Galilean invariance, as well as an associated momentum conservation law, which leads to a number of important tools, such as the Morawetz and virial identities, which do not have
direct analogues for gKdV\footnote{The gKdV, on the other hand, enjoys the property that radiation can only move in a leftward direction, while solitons (in the focusing case) move to the right.  This fact is extremely useful in the stability theory of solitons for gKdV, but is not so useful for Conjecture \ref{gkdv-conj}, as can already be seen from the proof of Theorem \ref{main1} below.}.  

Our first main result is to partly explain this discrepancy by showing that Conjecture \ref{gkdv-conj} is in fact logically stronger (except for a factor of $\sqrt{\frac{6}{5}}$) than Conjecture \ref{nls-conj}.

\begin{theorem}[gKdV scattering implies NLS scattering]\label{main1} Fix $p=5$ and $\mu=\pm 1$, and suppose that Conjecture \ref{gkdv-conj} is true for masses $M(u) < m$ for some $m \in (0,+\infty)$.  Then Conjecture \ref{nls-conj} is also true for masses $M(u) < \sqrt{\frac{5}{6}} m$.
\end{theorem}

The proof of this theorem is in fact rather simple, following from the observation in \cite{cct} (see also \cite{BoydChen}, \cite{schneider}) that a certain modulated,
rescaled version of NLS solutions are approximately gKdV solutions.  The main difficulty is to justify the approximation procedure, though this is not too difficult after using the well-posedness theory.  We give this proof in Section \ref{approx-sec}.

The loss of $\sqrt{\frac{6}{5}}$ is not serious in the defocusing case, as it still allows us to see that the full defocusing version of Conjecture \ref{gkdv-conj} logically implies (and thus is at least as hard as) the 
full defocusing version of Conjecture \ref{nls-conj}.  However it is somewhat
inefficient in the focusing case, as it only lets us show that the full focusing version of Conjecture \ref{gkdv-conj}
is at least as hard as Conjecture \ref{nls-conj} in the small mass case $M(u) < \sqrt{\frac{5}{6}} M(Q)$.  This type of loss is
explained by the fact that our asymptotic embedding of NLS into gKdV does not map NLS solitons to gKdV solitons (which should not be expected in any event, given that the perturbation theories for the two soliton families are significantly different).

\subsection{Monotonicity formula}

We now return to the general gKdV equation (in which $p$ is not necessarily $5$), and inspect the conserved mass \eqref{mass}
and energy \eqref{energy} more closely.  If we define the \emph{mass density}
$$\rho(t,x) := u(t,x)^2$$
the \emph{mass current}
$$j(t,x) := 3 u_x(t,x)^2 + \frac{2\mu p}{p+1} |u(t,x)|^{p+1}$$
the \emph{energy density}
$$e(t,x) := \frac{1}{2} u_x(t,x)^2 + \frac{\mu}{p+1} |u(t,x)|^{p+1}$$
and the \emph{energy current}
$$k(t,x) := \frac{3}{2} u_{xx}(t,x)^2 + 2 \mu |u(t,x)|^{p-1} u_x(t,x)^2 + \frac{\mu^2}{2} |u(t,x)|^{2p}$$
then a routine computation verifies (for Schwartz solutions, at least) the pointwise conservation laws
\begin{align}
\rho_t + \rho_{xxx} &= j_x \label{m-law}\\
e_t + e_{xxx} &= k_x \label{e-law}.
\end{align}

Let us now define the \emph{centre-of-mass}
$$ \langle x \rangle_M := \frac{1}{M(u)} \int_\R x \rho(t,x)\ dx$$
and the \emph{centre-of-energy}
$$ \langle x \rangle_E := \frac{1}{E(u)} \int_\R x e(t,x)\ dx$$
assuming that the solution has non-zero mass and energy and has enough decay that these quantities are well-defined. We conclude from \eqref{m-law}, \eqref{e-law} and an integration by parts the equations of motion
\begin{equation}\label{pitm}
\partial_t \langle x \rangle_M = - \frac{1}{M(u)} \int_\R j(t,x)\ dx 
\end{equation}
and
\begin{equation}\label{pite}
 \partial_t \langle x \rangle_E = - \frac{1}{E(u)} \int_\R k(t,x)\ dx.
\end{equation}
In the defocusing cases, $j$ and $k$ are non-negative, and so we obtain the well-known fact that the centre-of-mass
and centre-of-energy always move to the left.  In contrast, in the focusing equations, the soliton solutions defined earlier show that the centre-of-mass\footnote{The centre-of-energy of a soliton also moves to the right in the sub-critical cases
$p < 5$.  In the critical case $p=5$ the energy is zero, and so the centre-of-energy concept does not make sense; in the supercritical cases the energy is negative and so the concept is probably useless.} can move to the right.

\begin{remark} In the case $p=3$ (i.e. the modified KdV equation), or when $\mu=0$, 
the mass current and energy density are related by the simple relationship $j = 6e$.  This shows that the centre of mass 
evolves linearly:
$$ \partial_t \langle x \rangle_M = - \frac{6 E(u)}{M(u)}.$$
In more general defocusing cases, we can only say that $\partial_t \langle x \rangle_M$ is negative and comparable in magnitude to $\frac{E(u)}{M(u)}$.  
\end{remark}

The second main result of this paper is that in the defocusing case, the centre-of-energy moves to the left faster than the centre-of-mass:

\begin{theorem}[Monotonicity formula]\label{main} Let $p \geq \sqrt{3}$ and $\mu = +1$.  Let $u: \R \times \R \to \R$ be a global Schwartz solution of \eqref{eq:gkdv} which is not identically zero.  Then we have
\begin{equation}\label{mono}
 \partial_t \langle x \rangle_M - \partial_t \langle x \rangle_E > 0.
\end{equation}
In particular, the quantity $\langle x \rangle_E - \langle x \rangle_M$ is strictly monotone decreasing in time.  Related to this, we have the dispersion estimate
\begin{equation}\label{disp}
 \sup_{t \in \R} \int_\R |x-x(t)| [\rho(t,x) + e(t,x)]\ dx = +\infty
\end{equation}
for any function $x: \R \to \R$.
\end{theorem}

\begin{remark} The hypothesis that $u$ is Schwartz can be relaxed somewhat, for instance to the weighted Sobolev space $H^{2,2}_x$.  This would be particularly relevant in the cases where $p$ is not an odd integer, as one does not expect Schwartz solutions in general in such cases.  
\end{remark}

\begin{remark} This theorem can be explained intuitively as follows.  From stationary phase one expects that portions of the solution $u$ at frequency $\xi$ should move left with velocity $-3\xi^2$, thus high frequency components should move left extremely fast compared to low frequency components.  The energy density is more weighted towards high frequencies than the mass density, which then explains the monotonicity.
\end{remark}

\begin{remark} The condition $p \geq \sqrt{3}$ is an artefact of our method; we do not have any particular explanation for this exponent and it may be possible to obtain a similar dispersion estimate for lower $p$.  However, this hypothesis already covers the most important cases $p=2,3,4,5$.
\end{remark}

This theorem rules out any strongly localised ``pseudosoliton'' solutions to the defocusing gKdV equation, even if they
move to the left; it is a very weak analogue of the ``Liouville theorem'' in \cite{liouville}, which was for the focusing case in the vicinity of a soliton, or to the result of \cite{martel} which holds for arbitrary large data in focusing or defocusing cases, (but with the concentration $x(t)$ assumed to be moving to the \emph{right}).  This result may thus be viewed as analogous to the Morawetz or virial inequalities for the defocusing nonlinear Schr\"odinger (NLS) equation.  In analogy with the theory of the critical NLS (see e.g. \cite{borg:book}, \cite{gopher}, \cite{rv}, \cite{visan}, \cite{tao-visan}, \cite{kenig}) it is thus likely that suitably localised versions of this monotonicity formula could be useful in establishing large data global decay estimates for these defocusing equations, for instance in the critical case $p=5$, although it appears that a substantial amount of 
additional work would be needed to accomplish this.  

We thank Frank Merle for helpful discussions, and Shuanglin Shao and Monica Visan for a correction.

\section{Notation and basic estimates}

We use the notation $X \lesssim Y$ to denote a bound of the form $X \leq CY$ for some absolute constant $C$.
We use the usual mixed spacetime Lebesgue norms
$$ \| u \|_{L^q_t L^r_x(I \times \R)} := (\int_I (\int_\R |u(t,x)|^r\ dx)^{q/r} dt)^{1/q}$$
and
$$ \| u \|_{L^r_x L^q_t(I \times \R)} := (\int_\R (\int_I |u(t,x)|^q\ dt)^{r/q} dx)^{1/r}$$
with the usual modifications when $q,r = \infty$.  We caution that $L^q_t L^r_x$ and $L^r_x L^q_t$ are not equal except when $q=r$; in such cases we abbreviate $L^q_t L^q_x = L^q_x L^q_t = L^q_{t,x}$.  We use $C^0_t L^2_x$ to denote the closed subspace
of $L^\infty_t L^2_x$ (with the same norm) which is continuous as a function from time to $L^2_x$.

We abbreviate $\partial_x u$ as $u_x$, $\partial_{xx} u$ as $u_x$, etc.  We also define the usual fractional derivative operators $|\nabla_x|^\alpha$ by requiring
$$ |\nabla_x|^\alpha u(x) = \int_\R |\xi|^\alpha \hat u(\xi) e^{ix\xi}\ d\xi$$
whenever $u$ is Schwartz and has the Fourier representation
$$ u(x) = \int_\R \hat u(\xi) e^{ix\xi}\ d\xi.$$

We need the Kato smoothing, Strichartz, and maximal function estimates
\begin{equation}\label{kato}
\begin{split}
&\| |\nabla_x|^{-1/4} u \|_{L^4_x L^\infty_t(I \times \R)} 
+ \| u \|_{L^5_x L^{10}_t(I \times \R)} + \| |\nabla_x|^{1/6} u \|_{L^6_{t,x}(I \times \R)}\\
&+ \| |\nabla_x| u \|_{L^\infty_x L^2_t(I \times \R)}
+ \| u \|_{C^0_t L^2_x(I \times \R)}
\\
&\quad \lesssim \| u(t_0) \|_{L^2_x(\R)} + \| F \|_{L^1_x L^2_t(I \times \R)}
+ \|G\|_{L^1_t L^2_x(I \times \R)}
\end{split}
\end{equation}
whenever $u,F$ are smooth solutions on $I \times \R$ to the equation
\begin{equation}\label{duh}
 \partial_t u + \partial_{xxx} u = \partial_x F + G
\end{equation}
for some time interval $I$, and $t_0$ is any time in $I$; see \cite{kpv}.  Observe that the second and third estimates in \eqref{kato} can be obtained by complex interpolation\footnote{The complex differentiation operators $|\nabla_x|^{it}$ do not preserve all of the left-hand norms.  However, they preserve the right-hand norms, and by commuting these operators with the equation \eqref{duh} we can justify the complex interpolation.  We omit the details.} from the first and fourth estimates, which are the maximal function estimate and sharp Kato smoothing estimate respectively.

\section{Embedding NLS inside gKdV}\label{approx-sec}

We now prove Theorem \ref{main1}.  
Fix $\mu$, and suppose that $m > 0$ is such that Conjecture \ref{gkdv-conj} holds for masses $M(u) < m$.
Let $u: I \times \R \to \C$ be a fixed Schwartz solution to the NLS equation \eqref{nls} of mass $M(u) < \sqrt{\frac{5}{6}} m$ on a compact interval $I$; without loss of generality we shall assume the time interval $I$ to contain $0$.  Let $N \gg 1$ be an extremely large parameter (eventually it will go to infinity).  We shall adopt the following asymptotic conventions:

\begin{itemize}
\item (Uniform bound) The notation $O(X)$ shall refer to an expression bounded in magnitude by $C( M(u), m ) X$, where $0 < C(M(u),m) < \infty$ is allowed to depend on the mass threshold $m$ and the mass $M(u)$, but not on the time interval $I$ or the solution $u$ itself.
\item (Non-uniform bound) The notation $O_{u,I}(X)$ shall refer to an expression bounded in magnitude by $C(m,u,I)X$, where $0 < C(m,u,I) < \infty$ is now also allowed to depend on the time interval $I$ and on the Schwartz solution $u$ (in particular, it can depend on various Schwartz seminorms of $u$).  Thus every expression of the form $O(X)$ is also of the form $O_{u,I}(X)$, but not conversely.
\item (Non-uniform decay) The notation $o(1)$ shall refer to an expression bounded in magnitude by $c(N, u, I, m)$, where we have $\lim_{N \to \infty} c(N,u,I,m) = 0$ for fixed $u,I,m$.  Thus $o(1)$ expressions decay to zero as $N \to \infty$, but the rate of decay may depend on the solution $I$.  Note that the product of an expression which is $O_{u,I}(1)$ and an expression which is $o(1)$ will continue to be $o(1)$.  Also, any expression which is $o(1)$ will also be $O(1)$ when $N$ is sufficiently large depending on $u,I,m$.
\end{itemize}

Using this asymptotic notation, our task is thus to show that
\begin{equation}\label{goal}
\|u\|_{L^6_{t,x}(I \times \R)} = O(1).
\end{equation}

The key algebraic observation (as in\footnote{The numerology is very slightly different in \cite{cct}, as one was considering the $p=3$ cases of both gKdV and NLS, but the idea is otherwise identical.} \cite{cct}) is as follows.  
We introduce the function $u_N: I \times \R \to \R$ defined by
$$ u_N(t,x) := 
\frac{8^{1/4}}{5^{1/4} N^{1/4}} \Re \left[ e^{iNx} e^{iN^3 t} u( t, \frac{x+3N^2 t}{3^{1/2} N^{1/2}} ) \right]$$
and observe that $u_N$ approximately solves the gKdV equation \eqref{eq:gkdv}.  Indeed, a rather tedious computation using \eqref{nls} shows that
$$
(\partial_t + \partial_{xxx}) u_N = \mu \partial_x(|u_N|^4 u_N) + E_1 + E_2 + E_3
$$
where the errors $E_j = E_{j,N}$ for $j=1,2,3$ are given by
\begin{align*}
E_1 &:= 
N^{-1/4} \sum_{k=3,5} C_{1,k} \Re \left[ e^{ikNx} e^{ikN^3 t} |u|^4 u( t, \frac{x+3N^2 t}{3^{1/2} N^{1/2}} ) \right] \\
E_2 &:= N^{-5/4} \sum_{k=1,3,5} C_{2,k} \Re \left[ e^{ikNx} e^{ikN^3 t} (|u|^4 u)_x( t, \frac{x+3N^2 t}{3^{1/2} N^{1/2}} ) \right]\\ 
E_3 &:= C_3 N^{-7/4} \Re \left[ e^{iNx} e^{iN^3 t} u_{xxx}( t, \frac{x+3N^2 t}{3^{1/2} N^{1/2}} ) \right]
\end{align*}
and $C_{1,3}, C_{1,5}, C_{2,1}, C_{2,3}, C_{2,5}, C_3$ are absolute constants (depending only on $\mu$) whose exact value will not be important for us.
The point is that the dangerously large $N^{11/4}$ and $N^{5/4}$ terms have been completely cancelled, and the ``resonant'' term $j=1$ of the $N^{-1/4}$ term has also been cancelled (this is where the strange $\frac{8^{1/4}}{5^{1/4}}$ factor arises).

A rather direct computation (using the compactness of $I$ and the Schwartz nature of $u$) shows the qualitative bound
\begin{equation}\label{un}
\| u_N \|_{L^5_x L^{10}_t(I \times \R)} = O_{u,I}(1)
\end{equation}
(the point being that the right-hand side is uniform in $N$).
Indeed, $u_N$ is of size $O_{u,I}(N^{-1/4})$ on the region $\{ (t,x): t \in I, x = - 3N^2 t + O( N^{1/2} ) \}$, and is rapidly decaying away from this region.

The mass of $u_N(0)$ can be computed by a change of variables as
$$ \int_\R u_N(0,x)^2\ dx = \sqrt{\frac{6}{5}} \int_\R |u(0,x)|^2 + \Re( e^{2i 3^{1/2} N^{3/2} x} u(0,x)^2 )\ dx.$$
By the Riemann-Lebesgue lemma we thus have
$$ \int_\R u_N(0,x)^2\ dx = \sqrt{\frac{6}{5}} M(u) + o(1) < m$$
when $N$ is sufficiently large depending on $u,I,m$.  
Thus, by hypothesis (and the gKdV well-posedness theory, see \cite{kpv}), we may (for $N$ sufficiently large) find a solution $\tilde u_N: I \times \R \to \R$ to the gKdV equation \eqref{eq:gkdv} with initial data 
$\tilde u_N(0) = u_N(0)$, and with the uniform spacetime bound
\begin{equation}\label{tun}
 \| \tilde u_N \|_{L^5_x L^{10}_t(I \times \R)} = O(1).
\end{equation}
Also we have the mass bound
$$ M(\tilde u_N) = \int_\R u_N(0,x)^2 = O(1).$$
From \eqref{kato} and \eqref{eq:gkdv} we thus have
$$  \| |\nabla_x|^{1/6} \tilde u_N \|_{L^6_{t,x}(I \times \R)}
\lesssim O(1) + \| |\tilde u_N|^4 \tilde u_N \|_{L^1_x L^2_t(I \times \R)}$$
and hence by H\"older's inequality and \eqref{tun}
\begin{equation}\label{tun6}
\| |\nabla_x|^{1/6} \tilde u_N \|_{L^6_{t,x}(I \times \R)} = O(1).
\end{equation}

The next step is to pass from control of $\tilde u_N$ to control of $u_N$.  To do this we must first get some better
control on the error terms $E$.  We need some error estimates:

\begin{lemma} Let $e: I \times \R \to \R$ solve the forced Airy equation $\partial_t e + \partial_{xxx} e = E_1 + E_2 + E_3$ with
$e(0,x) = 0$.  Then we have
\begin{equation}\label{eun}
\| |\nabla_x|^{1/6} e \|_{L^6_{t,x}(J \times \R)} + \| e \|_{L^5_x L^{10}_t(I \times \R)} = o(1).
\end{equation}
\end{lemma}

\begin{proof} By linearity and the triangle inequality
it suffices to prove this claim for $E_1,E_2,E_3$ separately.  To prove the claims for $E_2,E_3$, observe from 
\eqref{kato} that it will suffice to show that
$$ \| E_2 \|_{L^1_t L^2_x} + \|E_3\|_{L^1_t L^2_x} = o(1),$$
but this is easily verified (recall that $I$ is compact and $u$ is Schwartz).  The same argument does not quite work for $E_1$.  However, we can almost solve for $e$ in this case; if we make the ansatz
$$ e = e_1 + e_2$$
where
$$ e_1 := N^{-1/4} \sum_{k=3,5} 
C_{1,k} \Re \left[ \frac{e^{ikN^3 t} - e^{ik^3 N^3 t}}{(ikN)^3 + (ikN^3)}
e^{ikNx} |u|^4 u( t, \frac{x+3N^2 t}{3^{1/2} N^{1/2}} ) \right] 
$$
(note that as $k \neq 1$, the denominator will not vanish), then one easily verifies (from the compactness of $I$ and the Schwartz nature of $u$) that $e_2(0) = 0$ and
$$ \| (\partial_t + \partial_{xxx}) e_2 \|_{L^1_t L^2_x} = o(1)$$
and hence by \eqref{kato}
$$ \| |\nabla_x|^{1/6} e_2 \|_{L^6_{t,x}(J \times \R)}  + \| e_2 \|_{L^5_x L^{10}_t(I \times \R)} = o(1).$$
Also, direct computation also shows
$$ \| |\nabla_x|^{1/6} e_1 \|_{L^6_{t,x}(J \times \R)} +
\| e_1 \|_{L^5_x L^{10}_t(I \times \R)} = o(1)$$
(basically because the denominator is as large as $N^3$) and the claim follows.
\end{proof}

If we now make the ansatz
$$ u_N = \tilde u_N + v + e$$
then $v(0) = 0$ and from construction of $e, \tilde u_N$ we have
$$ \partial_t v + \partial_{xxx} v = \partial_x( (\tilde u_N+v+e)^5 - \tilde u_N^5 ).$$
From \eqref{un}, \eqref{eun}, \eqref{tun} we have the qualitative bound
\begin{equation}\label{vun}
\| v\|_{L^5_t L^{10}_x(I \times \R)} = O_{u,I}(1).
\end{equation}
Actually we can improve this substantially.  Applying \eqref{kato}, we see that 
for any time interval $J \subset I$ and time $t_J \in J$, we have
\begin{align*}
\| v \|_{C^0_t L^2_x(J \times \R)} &+
\| |\nabla_x|^{1/6} v \|_{L^6_{t,x}(J \times \R)} +
\| v \|_{L^5_x L^{10}_t(J \times \R)}\\ 
&\lesssim \|v(t_J)\|_{L^2_x(\R)} + \| (\tilde u_N+v+e)^5 - \tilde u_N^5 \|_{L^1_x L^2_t(J \times \R)}.
\end{align*}
Now from the pointwise estimate
$$ (\tilde u_N+v+e)^5 - \tilde u_N^5 = O( |v|^5 ) + O( |e|^5 ) + O( |v| |\tilde u_N|^4 ) + 
O( |e| |\tilde u_N|^4 ) $$
and \eqref{tun}, \eqref{eun}, \eqref{vun} and H\"older we have
$$
\| (\tilde u_N+v+e)^5 - \tilde u_N^5 \|_{L^1_x L^2_t(J \times \R)}
\lesssim \| v \|_{L^5_x L^{10}_t(J \times \R)} \|\tilde u_N\|_{L^5_x L^{10}_t(J \times \R)} + \|v\|_{L^5_x L^{10}_t(J \times \R)}^5 + o(1)$$
and hence we have the inequality
\begin{align*}
\| v \|_{C^0_t L^2_x(J \times \R)} &+
\| |\nabla_x|^{1/6} v \|_{L^6_{t,x}(J \times \R)} +
\| v \|_{L^5_x L^{10}_t(J \times \R)}\\ 
&\lesssim \|v(t_J)\|_{L^2_x(\R)} + \| v \|_{L^5_x L^{10}_t(J \times \R)} \|\tilde u_N\|_{L^5_x L^{10}_t(J \times \R)} + \|v\|_{L^5_x L^{10}_t(J \times \R)}^5 + o(1).
\end{align*}
Standard continuity arguments thus show that $\|v(t_J)\|_{L^2_x(\R)} = o(1)$ and
 $\|\tilde u_N\|_{L^5_x L^{10}_t(J \times \R)}$ is smaller than a suitable absolute constant $\epsilon_0 > 0$, and
 $N$ is sufficiently large, then
$$ \| v \|_{C^0_t L^2_x(J \times \R)} + 
\| |\nabla_x|^{1/6} v \|_{L^6_{t,x}(J \times \R)} + \| v \|_{L^5_x L^{10}_t(J \times \R)} = o(1).$$
If we then use \eqref{tun} to subdivide $I$ into $O(1)$ intervals $J$ 
on which the $L^5_x L^{10}_t$ norm of $\tilde u_N$ is indeed less than $\epsilon_0$, then an easy induction followed by the triangle inequality thus ensures that
$$ \| v \|_{C^0_t L^2_x(I \times \R)} + \| |\nabla_x|^{1/6} v \|_{L^6_{t,x}(I \times \R)} + 
\| v \|_{L^5_x L^{10}_t(I \times \R)} = o(1).$$
when $N$ is sufficiently large.  Combining this with \eqref{tun6}, \eqref{eun} and the triangle inequality we obtain in particular a quantitative bound for $u_N$ when $N$ is sufficiently large:
$$ \| |\nabla_x|^{1/6} u_N \|_{L^6_{t,x}(I \times \R)} = O(1).$$
By definition of $u_N$, we thus have
$$
\| \Re e^{iN^3 t} |\nabla_x|^{1/6} \left[ e^{iNx} u( t, \frac{x+3N^2 t}{3^{1/2} N^{1/2}} ) \right] \|_{L^6_{t,x}(I \times \R)} 
= O( N^{1/4} ).$$
Making the change of variables $y := \frac{x+3N^2 t}{3^{1/2} N^{1/2}}$, this becomes
\begin{equation}\label{esp}
\| \Re e^{-2iN^3 t} |\nabla_y|^{1/6} \left[ e^{3^{1/2}iN^{3/2} y} u( t, y ) \right] \|_{L^6_{t,y}(I \times \R)} 
= O( N^{1/4} ).
\end{equation}
Now consider the commutator
$$ |\nabla_y|^{1/6} \left[ e^{3^{1/2}iN^{3/2} y} u( t, y ) \right] - 
[3^{1/2} N^{3/2}]^{1/6} e^{3^{1/2}iN^{3/2} y} u(t,y).$$
This has a Fourier expansion
$$ \int_\R e^{iy (\xi+3^{1/2} N^{3/2})} [ |\xi + 3^{1/2} N^{3/2}|^{1/6} - (3^{1/2} N^{3/2})^{1/6} ] \hat u(t,\xi)\ d\xi.$$
Some elementary stationary phase (using the Schwartz nature of $u$) shows that this quantity is $o(1)$ and is also rapidly decreasing in $y$, thus the net
contribution of this commutator to \eqref{esp} is $o(1)$.  Thus we have
$$
\| \Re e^{-2iN^3 t} [3^{1/2} N^{3/2}]^{1/6}  \left[ e^{3^{1/2}iN^{3/2} y} u( t, y ) \right] \|_{L^6_{t,y}(I \times \R)} 
= O( N^{1/4} )$$
and hence 
$$
\| \Re \left[  e^{-2iN^3 t} e^{3^{1/2}iN^{3/2} y} u( t, y ) \right] \|_{L^6_{t,y}(I \times \R)}^6 
= O( 1 ).$$
Expanding out the left-hand side and taking limits as $N \to \infty$ using the Riemann-Lebesgue lemma, we conclude
\eqref{goal} as desired.  This proves Theorem \ref{main1}.

\section{Proof of monotonicity formula}

In this section we prove Theorem \ref{main}.  We begin with the simpler task of proving \eqref{mono}. 
We use \eqref{pitm}, \eqref{pite} and reduce to showing that
$$
M(u) \int_\R k(t,x)\ dx > E(u) \int_\R j(t,x)\ dx$$
for each time $t$.  

Fix $t$. We introduce the strictly positive quantities $a,b,q,r,s$ by solving the equations
\begin{align*}
a^2 M(u) &:= \int_\R u_{xx}(t,x)^2\ dx \\
b^2 M(u)  &:= \int_\R |u(t,x)|^{2p}\ dx \\
aq M(u) &:= \int_\R u_x(t,x)^2\ dx \\
br M(u) &:= \int_\R |u(t,x)|^{p+1}\ dx \\
abs M(u) &:= p \int_\R |u(t,x)|^{p-1} u_x(t,x)^2\ dx.
\end{align*}
Our task can then be rewritten as
\begin{equation}\label{alg}
 \frac{3}{2} a^2 + 2 abs + \frac{1}{2} b^2 > (\frac{1}{2} aq + \frac{1}{p+1} br) (3aq + \frac{2p}{p+1} br).
\end{equation}
To prove this we observe a basic property of the quantities $q,r,s$ (which will help explain our notational conventions):

\begin{lemma}  The real symmetric matrix
$$
\begin{pmatrix}
1 & q & r \\
q & 1 & s \\
r & s & 1
\end{pmatrix}
$$
is positive semi-definite.
\end{lemma}

\begin{proof} For any real numbers $\alpha,\beta,\gamma$ we clearly have
$$ \int_\R | \gamma u(t,x) - \frac{\alpha}{a} u_{xx}(t,x) + \frac{\beta}{b} |u(t,x)|^{p-1} u(t,x) |^2\ dx \geq 0.$$ 
The left-hand side can be expanded out using integration by parts as
$$ M(u) [ \gamma^2 + \alpha^2 + \beta^2 - 2 \gamma \alpha r - 2 \gamma \beta s - 2 \alpha \beta t ]$$
and the claim follows.
\end{proof}

Taking determinants and minors we thus conclude that
\begin{equation}\label{rst}
0 < q,r,s \leq 1 \hbox{ and } 1 - q^2 - r^2 - s^2 + 2qrs \geq 0.
\end{equation}

We now are left with the purely algebraic task of deducing \eqref{alg} from \eqref{rst}.  
Expanding out, we reduce to showing that
\begin{equation}\label{alg2}
 \frac{3}{2} a^2 (1 - q^2) + ab (2s - \frac{p+3}{p+1} qr) + \frac{1}{2} b^2 (1 -  \frac{4p}{(p+1)^2} r^2) > 0.
 \end{equation}
Since $p > 1$, we have $\frac{4p}{(p+1)^2} < 1$; since $r$ is strictly positive, we may reduce to showing that
\begin{equation}\label{expand} 
 \frac{3}{2} a^2 (1 - q^2) + ab (2s - \frac{p+3}{p+1} qr) + \frac{1}{2} b^2 (1 -  r^2) \geq 0.
\end{equation}
By the quadratic formula (and the positivity of $a,b$), it will suffice to show that
$$ \frac{p+3}{p+1} qr - 2s \leq \sqrt{ 3 (1 - q^2) (1 - r^2) }.$$
Since $p \geq \sqrt{3}$, we have $\frac{p+3}{p+1} \leq \sqrt{3}$.  Thus it suffices to show that
$$ s \geq \frac{\sqrt{3}}{2} [qr - \sqrt{(1-q^2) (1-r^2)}].$$
Now from completing the square in the final inequality in \eqref{rst} we have
$$ (s-qr)^2 \leq (1-q^2) (1-r^2)$$
and thus
$$ s \geq qr - \sqrt{(1-q^2) (1-r^2)}.$$
On the other hand, $s$ is non-negative.  Since $\frac{\sqrt{3}}{2} \leq 1$, the claim follows.  This proves \eqref{mono}.

In order to prove the dispersion estimate \eqref{disp} we shall need to strengthen \eqref{mono} slightly.
Once again we shall need some asymptotic notation, this time a little different from in the preceding section. 

\begin{itemize}
\item We use $O(X)$ to denote any quantity bounded in magnitude by $C(u,p) X$, where $C(u)$ depends only on the solution $u$ and the exponent $p$.
\item We use $\Omega(X)$ to denote any positive quantity bounded \emph{below} by $C(u,p)^{-1} X$.
\item We use $\Theta(X)$ to denote a positive quantity which is both $O(X)$ and $\Omega(X)$.
\end{itemize}

Since $u$ is Schwartz and non-zero, we have
$$ M(u), E(u) = \Theta(1).$$
From the Gagliardo-Nirenberg inequality we thus have
$$ \int_\R |u(t,x)|^2\ dx, \int_\R |u_x(t,x)|^2\ dx = \Theta(1).$$
By Sobolev embedding we then conclude the pointwise bound $u(t,x) = O(1)$, and more generally
$$ (\int_\R |u(t,x)|^q\ dx)^{1/q} = O(1)$$
for $2 \leq q \leq \infty$.  Now we may assume for sake of contradiction
that \eqref{disp} fails, and so
\begin{equation}\label{rhoe}
 \int_\R |x-x(t)| (\rho(t,x) + e(t,x))\ dx = O(1).
\end{equation}
Since $\int_\R |u(t,x)|^2\ dx = \Theta(1)$, we conclude that
$$ \int_{x = x(t) + O(1)} |u(t,x)|^2\ dx = \Theta(1).$$
From the uniform bound on $u$ and H\"older's inequality we conclude that
$$ (\int_{x = x(t) + O(1)} |u(t,x)|^q\ dx)^{1/q} = \Theta(1)$$
for all $2 \leq q \leq \infty$, and thus
$$ (\int_\R |u(t,x)|^q\ dx)^{1/q} = \Theta(1).$$
By definitions of $b,r$ we now obtain the bounds
$$ br = \Theta(1)$$
(together with some additional bounds on $a,b,q,r,s$ which we will not need here).
We may then review the deduction of \eqref{alg2} from \eqref{expand} and
obtain the slight improvement
$$
 \frac{3}{2} a^2 (1 - q^2) + ab (2s - \frac{p+3}{p+1} qr) + \frac{1}{2} b^2 (1 -  \frac{4p}{(p+1)^2} r^2) > \Omega(1).$$
Working backwards (and recalling that $M(u) = \Theta(1)$ we conclude that
$$ M(u) \int_\R k(t,x)\ dx > E(u) \int_\R j(t,x)\ dx + \Omega(1)$$
and hence (since $M(u), E(u) = \Theta(1)$)
$$
 \partial_t \langle x \rangle_M - \partial_t \langle x \rangle_E > \Omega(1).$$
In particular we see that
$$\langle x \rangle_M - \langle x \rangle_E$$
cannot be bounded uniformly in time.  However, from \eqref{rhoe} (and the bounds $M(u), E(u) = \Theta(1)$) we easily
see that
$$ \langle x \rangle_M, \langle x \rangle_E = x(t) + O(1).$$
These two facts are contradictory, and we are done.

\begin{remark} A more careful evaluation of the above argument shows that we in fact have the more quantitative refinement
$$
\sup_{t \in I} \int_\R |x-x(t)| [\rho(t,x) + e(t,x)]\ dx \geq \Omega( |I|^c )$$
for all compact time intervals $I$ with $|I| \geq 1$, all functions $x: I \to \R$, and some explicitly computable
exponent $c = c(p) > 0$ depending only on $p$.  Furthermore, the implied constant in the $\Omega()$ notation depends only on the mass and energy of $u$ and on $p$.
We leave the details to the reader.
\end{remark}

\begin{remark} This phenomenon of separation of mass and energy does not seem to be particularly compatible with the embedding of NLS into gKdV used in Section \ref{approx-sec}.  The problem is that all the bounds in this section depend on the energy, however the solutions $\tilde u_N$ considered in Section \ref{approx-sec} have very large energy (comparable to $N^2$).  The energy density is in fact very close (in a weak sense) to $\frac{1}{2} N^2$ times the mass density, and so the centre-of-mass and centre-of-energy are virtually identical.
\end{remark}

\end{document}